\documentclass[epic,eepic,11pt]{amsart}

\usepackage{amsfonts,mathrsfs}
\usepackage[mathscr]{eucal}

\usepackage{amssymb}
\usepackage{amscd}
\usepackage{fancyhdr}
\usepackage[all,cmtip]{xy}

\usepackage{euler,eucal}

\DeclareMathAlphabet{\mathbf}{T1}{ppl}{bx}{n}
\DeclareMathAlphabet{\mathrm}{T1}{ppl}{m}{n}

\numberwithin{equation}{section}

\newcommand\note[1]%
{$^\dagger$\marginpar{\footnotesize{$^\dagger${#1}}}}

\def\({\left(}
\def\){\right)}
\def\<{\left<}
\def\>{\right>}

\newtheorem{theorem}{Theorem}[section]
\newtheorem{proposition}[theorem]{Proposition}
\newtheorem{lemma}[theorem]{Lemma}
\newtheorem{definition}[theorem]{Definition}

\newtheorem{corollary}[theorem]{Corollary}

\theoremstyle{definition}
\newtheorem{example}[theorem]{Example}
\newtheorem{remark}[theorem]{Remark}

\newcommand\lie{\mathfrak}

\newcommand\g{\lie{g}}

\newcommand\bb[1]{{\text{\bf#1}}}

\newcommand\Z{\bb{Z}}

\newcommand\R{\mathbb{R}}
\newcommand\C{\mathbb{C}}

\newcommand\CP{\mathbb{CP}}

\newcommand\J{\mathcal{J}}
\newcommand\cO{\mathcal{O}}

\newcommand     {\comment}[1]   {}
\newcommand{\mute}[2] {}
\newcommand     {\printname}[1] {}

\newcommand\func[1]{\operatorname{\mathrm{#1}}}
\newcommand\funclim[1]{\operatorname*{\mathrm{#1}}}

\renewcommand\dim{\func{dim}}

\renewcommand\lim{\funclim{lim}}
\renewcommand\log{\func{log}}

\newcommand\rank{\func{rank}}

\newcommand\sur{\mathrel{\to\kern-1.8ex\to}}
\newcommand\iso{\mathrel{\hookrightarrow\kern-1.8ex\to}}

\newcommand\longhookrightarrow{\lhook\joinrel\longrightarrow}

\newcommand\longsur{\mathrel{\longrightarrow\kern-1.8ex\to}}
\newcommand\longiso{\mathrel{\longhookrightarrow\kern-1.8ex\to}}

\begin{document}

\bibliographystyle{amsalpha}
\date{\today}
\title{Generalized complex hamiltonian torus actions: Examples and constraints.
}

\author{Thomas Baird, Yi Lin}
\maketitle

\begin{abstract}
Consider an effective Hamiltonian torus action $T\times M \rightarrow M$ on a topologically twisted,
generalized complex manifold $M$ of dimension $2n$. We prove that the $rank(T) \leq n-2$ and that the
 topological twisting survives Hamiltonian reduction. We then construct a large new class of such actions
  satisfying $rank(T) = n-2$, using a surgery procedure on toric manifolds.
\end{abstract}
\section{Introduction}

Generalized complex geometry, developed by Hitchin \cite{H02} and Gualtieri \cite{Gua03}, forms a common generalization of both symplectic and complex geometry and provides a natural geometric framework
  for the understanding of certain recent developments in string theory. A generalized complex structure on a manifold $M$ incorporates a "twist" by a closed differential 3-form, $H \in \Omega^3(M ; \R)$, $dH =0$. When $H$ is not exact, so that $[H] \neq 0 \in H^3(M;\R)$,
    then we say that $M$ is \emph{topologically twisted}.

In \cite{LT05}, Lin and Tolman introduced Hamiltonian actions for (twisted) generalized complex manifolds,
generalizing symplectic Hamiltonian actions (see also \cite{BCG05}, \cite{Hu05}, \cite{Va05}, \cite{SX05}.).
In this paper we introduce some new examples of, and constraints upon, compact, topologically twisted generalized
 complex Hamiltonian torus actions. Our first result is

\begin{theorem}\label{th1}
Let $M$ be a connected, compact topologically twisted generalized complex Hamiltonian $T$-space where $T$ is a torus acting effectively.
 If $\dim(M) = 2n$, then $\rank(T) \leq n-2$.
\end{theorem}

This contrasts with the (untwisted) symplectic setting, where effective Hamiltonian torus actions may have rank
 as large as half the dimension of the manifold.

A compact generalized Hamiltonian $T$-manifold $M$ consists of a compact torus action $T \times M \rightarrow M$,
on an $H$-twisted generalized complex manifold $M$ and a $T$-invariant smooth map $$ \mu : M \rightarrow \lie{t}^*$$
 satisfying some compatibility conditions. The compatibility condition relevant for us is as follows. The generalized complex
  structure induces naturally a bundle endomorphism $ L \in \Gamma( End( T^* M))$ and we insist that the equivariant $1$-form,
   $$\alpha = L( d\mu) \in \Omega^1(M)^T \otimes \lie{t}^*$$ determines an equivariantly closed extension of $H$ in the Cartan complex
 :$$H +\alpha \in \Omega_T^3(M), ~~~d_T(H+\alpha) = 0.$$ In particular, $H+\alpha$ determines an equivariant cohomology
     class $[H +\alpha] \in H^3_T(M)$.

For each value $\theta \in \lie{t}^*$ of the moment map $\mu$ such that $T$ acts freely on $\mu^{-1}(\theta)$,
 Lin-Tolman define a twisted generalized complex structure
 on the orbit space $$M //_{\theta} T := \mu^{-1}(\theta) / T$$ which is called the generalized complex quotient in analogy with symplectic
  terminology. The twisting $3$-form $\widetilde{H} \in \Omega( M //_{\theta} T)$ determines a cohomology class lying in the image of $[ H + \alpha]$
   under the map
$$ H_T(M) \rightarrow H_T(\mu^{-1}(\theta)) \cong H( M //_{\theta} T)$$
induced by inclusion. In \cite{Lin07}, Lin observed that if $H^{\text{odd}}(M)=0$, then necessarily,
 the twisting $\widetilde{H}$ on the quotient is exact. We first observe that this statement has the following strengthening under the additional assumption that $M$ is compact.

 \begin{proposition} \label{strengthening} Suppose that $M$ is a compact generalized complex Hamiltonian $T$-space with moment map $\mu$ and an exact twisting three form $H$, and that $T$ acts freely on the level set $\mu^{-1}(0)$. Then the generalized complex quotient $M//_0T$ has an exact twisting three form $\widetilde{H}$.
 \end{proposition}

 Our second main result is a converse of Propostion \ref{strengthening}

\begin{theorem}\label{th2}
For $T = U(1)$, let $M$ be a compact Hamiltonian $T$-manifold with topological twisting $H$ with moment map $\mu$ and let $M_{\theta} := {\mu}^{-1}(\theta)$. For any regular value $\theta \in \lie{t}^*$, the induced map $H_T(M) \rightarrow H_T(M_{\theta})
  \cong H(M//_{\theta} T)$ sends $[H + \alpha]$ to something nonzero. In particular, when $T$ acts freely on $M_{\theta}$,
  the generalized complex quotient $M//_{\theta} T$ is topologically twisted.
\end{theorem}

Using reduction in stages we deduce,

\begin{corollary}\label{ld}
Let $M$ be a Hamiltonian $T$-manifold with topological twisting $H$ and moment map $\mu$, and suppose that
the action of $T$ is \textbf{quasi-free}, i.e.,  all finite stabilizers $T_x \subset T$ are trivial.
For any regular value $\theta \in \lie{t}^*$, the generalized
complex quotient $M//_{\theta} T$ is topologically twisted.
\end{corollary}

The quasi-free hypothesis in Corollary \ref{ld} is needed to ensure that the intermediate stages
 of the reduction are actually smooth manifolds rather than mere orbifolds. We expect that Corollary \ref{ld}
  holds without this extra hypothesis.

The principal topological lemma that we make use of is the following result of Yasufumi Nitta \cite{NY07}, reformulated in the current form in \cite{BL08}.

\begin{proposition} \label{moment-map-Morse-Bott}
Consider the Hamiltonian action of a compact torus $T$ on a compact twisted GC-manifold $M$ with a generalized moment map
$\mu: M \rightarrow \lie{t}^*$. Then $\mu$ is a nondegenerate abstract moment map in the sense of \cite{GGK02}
(see Definition \ref{nondege-abstract-moment-map}).
\end{proposition}

Nondegenerate moment maps were introduced by \cite{GGK02} in order to axiomatize the Morse theoretic properties of symplectic moment maps.
 As noted in \cite{GGK02}, many important phenomena in symplectic geometry, such as Atiyah-Guillemin-Sternberg convexity, Kirwan surjectivity,
  and localization are often due to the properties of abstract moment maps more than to the symplectic structures. In view of this, it is not
   surprising that Proposition \ref{moment-map-Morse-Bott} has many geometric and topological consequences, some of which are explored in
    \cite{NY07}, \cite{BL08}.

Our strategy to prove Theorem \ref{th1} and Theorem \ref{th2} will be to apply the Morse theoretic techniques developed by Tolman and Weitsman
 in \cite{TW98} to generalized complex moment maps. Though their paper focused on the symplectic category,  it was noted in \cite{TW98} that the results generalize easily to nondegenerate abstract moment maps. We apply these methods to achieve different aims than those for which they were developed.

Since Theorems \ref{th1} and \ref{th2} concern topologically twisted Hamiltonian actions, it is natural to ask if interesting examples of such actions exist. In \cite{Lin07}, Lin constructed the first compact examples
    of topologically twisted generalized complex Hamiltonian manifolds which arise as symplectic fibrations over topologically twisted complex manifolds.

In \S \ref{examples}, we construct a new class of compact examples of topologically twisted generalized complex Hamiltonian torus actions.
 Our construction begins with a symplectic toric manifold $X$ of dimension $2n-2$ satisfying some conditions, and involves a surgery
  construction on the product $ T^2 \times X$, where $T^2$ is a 2-torus. In the process one of the $S^1$ symmetries
   is broken and we are left with a rank $n-2$ torus action realizing the upper bound of Theorem \ref{th1}. The surgery is an equivariant
    version of the surgery Gualtieri and Cavalcanti introduced in \cite{CaGu07} to produce topologically twisted generalized complex 4-manifolds from
     symplectic 4-manifolds.

\section{Review of generalized complex geometry}

Let $M$ be a manifold of dimension $n$. There is a natural pairing
of type $(n,n)$ which is defined on $TM\oplus T^*M$ by
\[ \langle X+\alpha, Y+\beta \rangle
=\dfrac{1}{2}\left(\beta(Y)+\alpha(X)\right),\] and which extends
naturally to $T_{\C}M\oplus T_{\C}^*M$.

For a closed three form $H$, the $H$-\textbf{twisted Courant
bracket} of $T_{\C}M\oplus T^*_{\C}M$ is defined by the identity
\[
[X+\xi,Y+\eta]=[X,Y]+L_X\eta-L_Y\xi-\dfrac{1}{2}d\left(\eta(X)-\xi(Y)\right)+\iota_Y\iota_XH.\]

      Let $\sigma$ be the linear map on $\wedge T^*M$ which acts on
      decomposables by
      \begin{equation}\label{anti-auto1} \sigma (f_1\wedge f_2 \wedge
      \cdots \wedge f_q)=f_q \wedge f_{q-1}\wedge \cdots \wedge
      f_1,\end{equation} then we have the following  pairing,
      called the Mukai pairing \cite{Gua07}, defined on differential forms:
      \[ (\xi_1,\xi_2)=(\sigma(\xi_1)\wedge \xi_2)_{\text{top}}, \] where
      $(\cdot,\cdot)_{\text{top}}$ indicates taking the top degree
      component of the form.

A \textbf{generalized almost complex structure} on a manifold $M$ is
an orthogonal bundle map $\mathcal{J}:TM\oplus T^*M \rightarrow
TM\oplus T^*M$ such that $\mathcal{J}^2=-1$. Moreover, $\mathcal{J}$
is an $H$-\textbf{twisted generalized complex structure} if the
sections of the $\sqrt{-1}$ eigenbundle of $\mathcal{J}$ is closed
under the $H$-twisted Courant bracket.

    \begin{example} (\cite{H02}, \cite{Gua07})\label{symplectic-g-calabi-Yau}
     Let $(M,\omega)$ be a  $2n$ dimensional symplectic manifold.
     Then

        \begin{equation} \label{symplectic}
        \mathcal{J}_{\omega}=\left(\begin{matrix} 0 & -\omega^{-1} \\\omega
        &0\\\end{matrix} \right): TM\oplus T^*M\rightarrow TM \oplus T^*M
        \end{equation}
        is a generalized complex structure on $M$.
\end{example}

  Let $B$ be a two-form on a manifold M, and consider the
  orthogonal bundle map defined by

  \[e^B=\left(\begin{matrix} 1 & 0\\ B&1 \end{matrix}\right): TM\oplus T^*M\rightarrow TM \oplus T^*M,\]
  where $B$ is regarded as a skew-symmetric map from $TM$ to $T^*M$ given by
  \[ X\longmapsto B(X,\cdot).\]

  If $\mathcal{J}$ is
  an $H$-twisted generalized complex structure on M, then
  $\mathcal{J}_B := e^B\mathcal{J} e^{-B}$ is another $H+dB$-twisted
  generalized complex structure on M, called the $B$-\textbf{transform} of
  $\mathcal{J}$.

 Locally, a generalized complex structure can be described by pure spinors.
Recall that the Clifford algebra of $C^{\infty}(TM\oplus T^*M)$ with the natural
pairing acts on differential forms by
\[ (X+\xi)\cdot \alpha=\iota_X\alpha+\xi\wedge \alpha.\]
A  differential form $\varphi $ on an open subset $U$ is called a pure
spinor if it is nowhere vanishing on $U$, and if for any $x\in U$,
\[L_{\varphi, x}=\{X+\xi\in T_{x}M\oplus T_x^*M, \, (X+\xi)\cdot \varphi=0      \}\]
 is a maximal isotropic subspace of $T_xM\oplus T_x^*M$.  $L_{\varphi}$ gives rise to
the $i$-eigenbundle of an unique $H$-twisted generalized complex structure $\J$
if and only if the pure spinor $\varphi$ satisfies the following two conditions
 \begin{itemize}\label{GC-spinor} \item [a)]    the Mukai pairing $(\varphi, \overline{\varphi})\neq 0$;
 \item [b)]
     $d\varphi- H\wedge \varphi= (X+\xi)\cdot \varphi $ for some $X+\xi \in C^{\infty}(TM\oplus T^*M)$.
      \end{itemize}

In this case, we also call $\varphi$ a pure spinor associated to the generalized complex structure $\J$. Conversely,
an $H$-twisted generalized complex structure $\J$ on a manifold $M$ determines a line bundle,
called the \textbf{canonical line bundle}, whose local non-zero sections are given by pure spinors satisfying
the above two conditions.

 The effect of $B$-transforms on the canonical line bundle of a
 generalized complex structure has been studied in \cite{Ca06}.

 \begin{lemma}(\cite{Ca06} )\label{Btransformoperator}
 Locally, let $\rho$ be a pure spinor associated to the $H$-twisted generalized complex structure $\mathcal{J}$, and
 $B$ a two form. Then $e^B\wedge \rho$ is the pure spinor associated to the $(H+dB)$-twisted generalized complex structure
 $\J_B$, the $B$-transform of $\J$.     \end{lemma}

    \begin{example} (\cite{H02})\label{symplectic-g-c-manifolds}
     Let $(M,\omega)$ be a  dimensional symplectic manifold, and $B$ a two form on $M$.
     Then $\varphi=e^{B+i\omega}$ is a pure spinor associated to the $B$-transform of the generalized complex
     structure  $\J_{\omega}$ as defined in (\ref{symplectic}), which is often called the $B$-transform of the
     symplectic structure $\omega$, or simply a $B$-symplectic structure.
   \end{example}

      At any point $x\in U$,  a pure spinor $ \varphi$ takes the form $e^{B+i\omega}\wedge \Omega$, where
       $B$ and $\omega$ are real $2$-forms, and $\Omega$ is a decomposable complex $k$-form. In general, the degree of
       $\Omega$ depends on $x$. If $\varphi$ is a pure spinor associated to a generalized complex structure
       $\J$, then the degree of $\Omega$ at $x$ is defined to be the \textbf{type} of the generalized complex structure $\J$ at $x$.
    The following example of a generalized complex structure with type changing is taken from
   \cite{CaGu07}, to which we refer for a detailed account.

\begin{example} \label{local-model}   Consider $\C^2$ with complex coordinates $z_1, z_2$. The
         differential form    \[\rho = z_1 + dz\wedge dz_2\]
         is a pure spinor and determines a generalized complex structure which
         has type $2$ along $z_1 = 0$ and type $0$ elsewhere.
         Note that $\rho$ is invariant under translations in the $z_2$ direction.
         Taking a quotient by the standard $Z^2$ action, we obtain a generalized
         complex structure $\J$ on the torus fibration $D\times T^2$, where $D$ is the unit open disc in the
         $z_1$-plane. While away from the central fibre $\{0\}\times T^2$, $\J$
         has type zero,  its type jumps from 0 to 2 along $\{0\}\times T^2$.

         More explicitly,  let
         \begin{equation} \label{B-symplectic}B = d \log r\wedge d\theta_2 -d\theta_1\wedge d\theta_3, \,\,\,\,\,
          \gamma =   d \log r \wedge d\theta_3 + d\theta_1
                \wedge d\theta_2, \end{equation}   where $r, \theta_1$ are the polar coordinates on $D$, and $\theta_2, \theta_3$
            are coordinates for $T^2$ with unit periods.  Then outside the type jumping locus, the generalized
             complex structure $\J$ is the $B$-transform of the symplectic structure $\gamma$.
\end{example}

\subsection{Review of equivariant cohomology and generalized complex quotients}

     Let $G$ be a compact connected Lie group, and $EG$ a contractible
     topological space on which $G$ acts freely. The equivariant cohomology (with complex coefficients) of
     the $G$-space $M$ is defined to be
     \[H_G(M):=H((EG\times M)/G; \C) .\]
     We would like to mention two fundamental properties of the equivariant cohomology theory.
     For a proof of them, we refer to \cite{GS99}.
     \begin{itemize}
    \item [a)]  Suppose that the action of $G$ on $M$ is locally free, i.e. $G$ acts on $M$ with
    finite stabilizers. Then
     \[ H_G(M)\cong  H(M/G).\]
     \item [b)] Suppose that $f: M \rightarrow N$ be an equivariant map between the $G$-spaces $M$ and $N$.
     Then $f$ induces a pullback homomorphism on the equivariant cohomology groups
     \[f^*:H_G(N)\rightarrow H_G(M).\]
      \end{itemize}

      Now suppose that $M$ is a $G$-manifold. Then its Cartan complex $\Omega_G(M)=(S\g^*\otimes \Omega(M))^G$
      is a graded complex with the total grading
      \[ \Omega_G^k(M)=\bigoplus_{2i+j=k}  (S^i\g^*\otimes \Omega^{j}(M))^G \]
      and the equivariant differential $d_G$ given by
       \begin{equation} \label{coboundary} (d_G\rho)(\xi)=d(\rho(\xi))-\iota_{\xi_M}\rho(\xi),
       \end{equation}
       where $\rho \in \Omega_G(M)$ is regarded as an equivariant polynomial from $\g$ to
        $\Omega(M)$, and $\iota_{\xi_M}$  denotes inner product with the vector field $\xi_M$ on $M$ induced
       by $\xi \in \g$.  The equivariant de Rham theorem (c.f. \cite{GS99}) asserts that the cohomology
         of the Cartan complex
         \[H(\Omega_G(M), d_G) \cong H_G(M,\C).\]


    We are ready to introduce the notion of Hamiltonian actions on generalized complex manifolds.
\begin{definition}
\cite{LT05}) \label{deftmm}
Let a compact Lie group $G$ with Lie algebra $\g$ act on a manifold
$M$, preserving an $H$-twisted generalized complex structure
$\mathcal{J}$, where $H \in \Omega^3(M)^G$ is closed.
The action of $G$ is said to be Hamiltonian if there exists a smooth
equivariant function $\mu:M \rightarrow \frak{g}^*$,  called the
{\bf generalized moment map}, and a $1$-form $\alpha \in
\Omega^1(M,\g^*)$, called the {\bf moment one form},  so that
\begin{itemize}
\item [a)]
$-\J d\mu^{\xi}=\xi_M+\alpha^{\xi}$ for all $\xi \in \g$, where $\xi_M$ denotes the
induced vector field.
\item [b)] $H+\alpha$ is an equivariantly closed three form in the
Cartan complex.
\end{itemize}
\end{definition}

Let a compact Lie group $G$ act on a twisted generalized complex
manifold $(M,\J)$ with generalized moment map $\mu$. Let $\cO_a$ be
the co-adjoint orbit through $a \in \g^*$. If $G$ acts freely on
$\mu^{-1}(\cO_a)$, then $\cO_a$ consists of regular values and $M_a
= \mu^{-1}(\cO_a)/G$ is a manifold, which is called the {\bf
generalized complex quotient}. The following result was
proved in \cite{LT05}.


\begin{proposition} \label{Twisted Complex Reduction}
Assume there is a Hamiltonian action of a compact Lie group $G$ on
an $H$-twisted generalized complex manifold $(M,\mathcal{J})$ with
generalized moment map $\mu: M \to \g^*$ and moment one-form $\alpha
\in \Omega^1(M,\g^*)$. Let $\cO_a$ be a co-adjoint orbit through $a
\in \g^*$ for which $G$ acts freely on $\mu^{-1}(\cO_a)$. Then
the generalized complex quotient
$M_a$ inherits an $\widetilde{H}$-twisted generalized complex
structure $\widetilde{\J}$; moreover, the cohomology class of $\widetilde{H}$
coincides with the image of $[H+\alpha] \in H_G(M)$ under the Kirwan map
\begin{equation} \label{kirwan-map} \begin{CD}
	\xymatrix{ H_G(M) \ar[r]^{i^*}  \ar[dr]^{k}    & H_G(\mu^{-1}(\cO_a)) \ar[d]^{\cong} \\
 & H(M_a)},
\end{CD} \end{equation}  where the horizontal map is the one induced by the inclusion
$ i: \mu^{-1}(\cO_a)\rightarrow M$. \end{proposition}

  \begin{example}\label{Hamiltonian-GCmflds}
  \begin{itemize}

  \item[a)]If $G$ acts on a symplectic manifold
  $(M, \omega)$ with moment map $\Phi$, then $G$ also preserves
 the generalized complex structure $\mathcal{J}_\omega$, and  $\Phi$ is a generalized moment map for
  this action.

 \item [b)] Let $G$ act on an $H$-twisted generalized complex manifold $(M,\J )$
  with generalized moment map $\mu$ and moment one-form $\alpha$.
  If $B \in \Omega^2(M)^G$, then $G$ acts on the $B$-transform of $(M, \J)$ with
  generalized moment map $\mu$ and moment one form $\alpha'$, where $(\alpha')^{\xi}=\alpha^{\xi}+\iota_{\xi_M}B$,
  $ \xi\in \frak{g}$.

  \item[c)] Suppose that $M_1$ and $M_2$ are two generalized complex manifolds,
  and that the action of a Lie group $G$ on $M_1$ is Hamiltonian. Then the action of
  $G$ on the generalized complex manifold $M_1\times M_2$ given by
  \[ g\cdot (x,y)=(g\cdot x, g\cdot y),\,\,\,\,\,\,(x,y)\in M_1\times M_2 \]
   is Hamiltonian.

  \end{itemize}
  \end{example}



\section{Review of nondegenerate abstract moment maps and Tolman-Weitsman techniques}

Following \cite{GGK02}, first we give a quick review of non-degenerate abstract moment maps. Let $G$ be a Lie group, $\frak{g}$ its Lie algebra,
 and $\frak{g}^*$ the dual space of $\frak{g}$. For any map $\Psi: M\rightarrow \frak{g}^*$ and any subgroup $H$ of $G$ with Lie algebra $\frak{h}$,
 we denote by $\Psi^H$ the composition of $\Psi$ with the natural projection $\frak{g}^* \rightarrow \frak{h}^*$. Similarly, for
 any Lie algebra element $\xi \in \frak{g}$ the $\xi$-th component of $\Psi$, i.e., the real function $<\Psi, \xi>$, is denoted by $\Psi^{\xi}$.

\begin{definition}\label{abstract-moment-map} Let $M$ be a $G$-manifold. An abstract moment map is a smooth map $\Psi: M\rightarrow \frak{g}^*$
 with the following properties
\begin{itemize}
\item [1.] $\Psi$ is $G$-equivariant;
\item [2.] For any subgroup $H$ of $G$, the map $\Psi^H: M\rightarrow \frak{h}^*$ is locally constant on the submanifold $M^H$ of points fixed by $H$.
\end{itemize}

\end{definition}

\begin{definition}\label{nondege-abstract-moment-map}  An abstract moment map $\Psi: M\rightarrow \frak{g}^*$ is non-degenerate if for every vector
 $\xi \in \frak{g}$,
\begin{itemize}

\item [1.] $\text{Crit}(\Psi^{\xi})=\{\xi_M=0\}$, and

\item [2.] $\Psi^{\xi}$ is a Morse-Bott function.

\end{itemize}

\end{definition}

For the remaining proofs, we will need to apply a theorem due to Tolman-Weitsman \cite{TW98}. Tolman and Weitsman proved the theorem
 for Hamiltonian circle actions on compact symplectic manifolds. But they also explained \cite[Remark 3.4]{TW98} that their argument
  works perfectly well for non-degenerate abstract moment maps.

\begin{theorem}[\cite{TW98}]\label{thtw}
Let $M$ be an $S^1$-manifold, and $f: M \rightarrow \R$ a non-degenerate abstract moment map for the $S^1$ action. Assume that $0$ is
 a regular value of $f$ and that $S^1$ acts with finite isotopy groups on $f^{-1}(0)$. Let $F$ denote the set of fixed points; write
 $M_-:= f^{-1}(-\infty,0)$, and $M_+:= f^{-1}(0,\infty)$. Define $$ K_{\pm}:=\{\alpha \in H^*_{S^1}(M)|\alpha|_{F\cap M_{\pm}} = 0\},$$
  and $$K := K_+ + K_- = K_+ \oplus K_-.$$
Then there is a short exact sequence: $$0\rightarrow K \rightarrow H_{S^1}^*(M)\xrightarrow{\kappa} H^*( M_{red}) \rightarrow 0,$$ where
$M_{red} := f^{-1}(0)/S^1$, and $\kappa$ is the Kirwan map as defined in (\ref{kirwan-map}).
\end{theorem}

\section{Constraint on the rank of the torus}

In this section we prove Theorem \ref{th1}. We will require the following result (see Corollary 5.9 of \cite{BL08}).

\begin{lemma}\label{afterthot}
Let $M$ be compact generalized Hamiltonian $T$-space, with moment 1-form $\alpha$. Then the restriction of $\alpha$ to the fixed point set $M^T$ is zero.
\end{lemma}

\begin{lemma}\label{doopy}
Let $M$ be a compact generalized Hamiltonian $T$-space with topological twist $H$. Then the
restriction of $H$ to the fixed point set $M^T$ is not exact.
\end{lemma}

\begin{proof}[Proof of Proposition \ref{strengthening} and Lemma \ref{doopy}]
It was observed in \cite{BL08} that the torus action of $T$ on $M$ is equivariantly formal in (untwisted) equivariant cohomology.
Thus the localization map $i^*:H_T^*(M) \rightarrow H_T^*(M^T)$ is injective, where $i: M^T \hookrightarrow M$ is inclusion of the
fixed point locus. Note that by Lemma \ref{afterthot}, the restricted equivariant de Rham form
$$ i^*(H + \alpha) = i^*(H) + i^*(\alpha)=i^*H \in \Omega(M^T)^T \otimes S(\lie{t}).$$
Moreover, because $T$ acts trivially on
$M^T$, we have $d_T = d$ in the Cartan model $\Omega_T(M^T)$, where we abuse notation by writing $d$ in place of $d \otimes id_{S\lie{t}^*}$.

If $H \in \Omega^3(M)$ is exact, it follows that $ i^*(H + \alpha)=i^*H$ is exact in $\Omega(M^T)^T \otimes S(\lie{t})$. So $H+\alpha$ is an exact equivariant three form in the Cartan model. By Proposition \ref{Twisted Complex Reduction}
the twisting three form $\widetilde{H}$ must be exact as well. This proves Proposition \ref{strengthening}.

 If $H\in\Omega^3(M)$ is non-exact, then $[H + \alpha] \neq 0$ in $H^3_T(M)$. This implies that the restricted equivariant de Rham form
$ i^*(H + \alpha) = i^*(H) \in \Omega(M^T)^T \otimes S(\lie{t})$ is not exact.  So we deduce that $i^*(H)$ is equivariantly closed and nonexact. Thus if there exists $\beta \in \Omega^2(M^T) = \Omega^2(M^T)^T$ satisfying
$d \beta = i^*(H)$, then $ d_T(\beta) = i^*(H) \in \Omega_T^3(M)$,
which is a contradiction. This completes the proof of Lemma \ref{doopy}.
\end{proof}

Because $M^T$ admits a nonexact 3-form, it follows that $H^3(M^T) \neq 0$ and consequently $M^T$ is at least 3 dimensional. Indeed we can do one better.

\begin{corollary}\label{doopy2}
Let $M$ be a compact $H$-twisted generalized Hamiltonian $T$-space and suppose that the twisting 3-form $H$ is not exact. Then some connected
component of $M^T$ has dimension at least four.
\end{corollary}

\begin{proof}
By Lemma \ref{doopy}, $M^T$ supports a nonzero 3-form so $M^T$ must be at least $3$-dimensional. According to \cite[Lemma 5.4]{Lin06},
$M^T$ is a generalized complex submanifold of $M$, so all components are even dimensional.
\end{proof}

Corollary \ref{doopy2} stands in stark contrast with the (untwisted) case of symplectic manifolds, where there is are abundant examples of Hamiltonian torus actions with isolated fixed points (e.g. toric varieties). Also striking is the conclusion
that all generalized Hamiltonian actions on (topologically twisted) compact generalized complex 4-manifolds must be trivial.
In fact, we get the stronger result:

\begin{proof}[Proof of Theorem \ref{th1}]
Let $F \subset M$ be a component of the fixed point set with dimension $4$ or greater. $F$ is an embedded, compact submanifold without boundary.
Let $\nu \rightarrow F$ be the normal bundle of $F$ in $M$ which we identify equivariantly with a tubular neighborhood of $F$. The $T$ action on
$\nu$ extends to $\nu \otimes \C$, which decomposes into a direct sum of eigenbundles $$ \nu \otimes \C = \oplus_{i=1}^m \nu_i \oplus \bar{\nu}_i$$
where $T$ acts on $\nu_i$ via a nonzero character $\chi_i: T \rightarrow U(1)$ and on the conjugate $\bar{\nu}_i$ via the conjugate character
$\bar{\chi}_i$. It follows that the subgroup of $K \subset T$ acting trivially on $\nu$ is equal to the kernel of the map:
$$ (\chi_1,...,\chi_m): T \rightarrow U(1)^m$$
Identifying $\nu$ with a neighborhood of $F$ in $M$, we find $\dim(M^K) = \dim(\nu) = \dim(M)$ and so $K$ acts trivially on $M$.
By hypothesis, $K$ is zero dimensional. Therefore:
$$ \dim(T) \leq m \leq (\dim(M) - \dim(F))/2 \leq n -2$$
\end{proof}

\section{topological twist of the quotient}

For a circle action, the Hamiltonian is simply a real valued function $f: M \rightarrow \R$. The sets $M_{max}$ and $M_{min}$, where $f$
achieves its maximum and minimum values respectively, are both components of the fixed point locus. In this section we refine our understanding
of Lemma \ref{doopy} to show that the restrictions of $H$ to both $M_{max}$ and $M_{min}$ are cohomologically nonzero. Theorem \ref{th2} then
follows from the Tolman-Weitsman description of the kernel of the Kirwan map.

\begin{lemma}\label{xmas1}
Let $T = U(1)$ be the circle group, and let $\pi: \nu \rightarrow M$ be a $T$-equivariant real vector bundle over a manifold. Consider the
commutative diagram

\begin{equation*}\begin{CD}
	\xymatrix{ M  \ar @/_/ [r]_{i} & \ar @/_/[l]_{\pi} \nu  & S(\nu) \ar[l]_{j}}
\end{CD}\end{equation*}
where $\pi$ is the bundle projection, $i$ is inclusion as the zero section and $j$ is an equivariant inclusion of the sphere bundle, $S(\nu)$.
Suppose that $\nu^T = i(M)$ and let $\eta_T \in \Omega_T(\nu)$ be a $d_T$-closed form such that $i^*(\eta_T)$ lies in
$ \Omega(M) \subset \Omega_T(M)$, i.e., $i^*\eta_T$ is an ordinary differential form.Then $\eta_T$ is exact if and only if $j^*(\eta_T)$ is exact.
\end{lemma}

\begin{proof}

Consider the commutative diagram of $T$-equivariant maps:

\begin{equation}\begin{CD}
\xymatrix{  ET \times M \ar[d]^{proj_M}& \ar[l]_{id\times \pi} ET \times \nu & \ar[l]_{id \times j}  ET \times S(\nu)\ar[d] \\
             M              & & S(\nu) \ar[ll]^{\pi \circ j}                 }
\end{CD}\end{equation}
This determines a map at the level of $T$-quotients and passing to cohomology determines a commutative diagram:

\begin{equation}\begin{CD}
\xymatrix{  H_T(M) \ar[r]^{(id \times \pi)^*}_{\cong} &  H_T(\nu) \ar[r] &   H_T(S(\nu)) \\
            H(M) \ar[u]^{proj_M^*}  \ar[rr]^{(\pi\circ j)^*}            & & H(S(\nu)/T) \ar[u]^{proj_{S(\nu)}^*}                }
\end{CD}\end{equation}
where the isomorphism $(id \times \pi)^*$ has inverse $(id \times i)^*$.

Now we know from the hypotheses that $[\eta_T] \in H_T(\nu)$ lies in the image of $(id \times \pi)^* \circ proj_M^*$. Thus it will suffice
to prove that $proj_{S(\nu)}^*$ and $(\pi \circ j)^*$ are injective.

Because $\nu^T = i(M)$, T acts with finite stabilizers on $S(\nu)$. It follows that $proj_{S(\nu)}^*: H(S(\nu)/T) \cong H_T(S(\nu))$
 is an isomorphism and is in particular injective.

Now choose $\nu_i \subset \nu$ to be a nontrivial eigenbundle of $\nu$ of weight $i$. Then it is possible to put a complex
structure $J$ on $\nu_i$ so that $T$ acts by scalar multiplication by a character of $T$. In particular, $S(\nu_i)/T $ is
homeomorphic to $P(\nu_i,J)$, the projective bundle associated to $(\nu_i,J)$.  We get maps:

\begin{equation}\begin{CD}
P(\nu_i,J) = S(\nu_i)/T @>>>  S(\nu)/T    @> \pi \circ j>>              M .
\end{CD}\end{equation}
The composed map is just projection. By the splitting principle (see \S 21 \cite{BT82}), the induced map
$H(M) \rightarrow H(P(\nu_i,J))$ is an injection, so $(\pi \circ j)^*$ must also be injective.

\end{proof}

\begin{lemma}\label{xmas2}
Let $M$ be compact Hamiltonian $T= U(1)$-manifold with topological twisting $H$, moment map $f$ and moment one form $\alpha$.
Let $F \subset M^T$ be a fixed point set component with $f(F)=c$. Then if the restriction of $H$ to $F$ is not exact in $\Omega^*(F)$,
then for sufficiently small $\epsilon >0$,  $H+\alpha$ restricts to a nonexact equivariant form on $M_{c+\epsilon}:=f^{-1}(c+\epsilon)$
and $M_{c-\epsilon}:=f^{-1}(c-\epsilon)$ when these are non-empty.
\end{lemma}

\begin{proof}
Suppose $c \neq f_{min}$, then the negative normal bundle $\nu_F \rightarrow F$ has positive rank.  We abuse notation by identifying $\nu_F$
 equivariantly with a subset of $M$. Consider the following commutative diagram:

\begin{equation}\begin{CD}
	\xymatrix{ H_T^*(M) \ar[rr] \ar[d] & & H_T^*(M_{c-\epsilon}) \ar[d]\\
H_T^*(F) \ar[r]^{\cong} & H_T^*(\nu_F) \ar[r]& H_T^*(\nu_F \cap M_{c-\epsilon})	\\
}
\end{CD}\end{equation}

For $\epsilon > 0$ sufficiently small, we have $M_{c-\epsilon} \cap \nu_F \cong S(\nu_F)$ a sphere bundle. Thus by Lemma \ref{xmas1},
$(H+\alpha)|_{M_{c-\epsilon} \cap \nu_F}$ is nonexact, so $\eta_T|_{M_{c-\epsilon}}$ is also nonexact. Similarly, for $M_{c+\epsilon}$.
\end{proof}

\begin{corollary}\label{snoopy}
Using hypotheses of Lemma \ref{xmas2}, let $M_{max}$ and $M_{min}$ be the sets on which $f$ achieves its minimum
and maximum respectively. Then $\eta$ restricts to a nonexact form on $M_{max}$ and $M_{min}$.
\end{corollary}

\begin{proof}
By Lemma \ref{doopy}, we know that $\eta_T$ restricts to a nonexact form on some component $F$ of $M^T$.
If $f(c) > f_{min}$ then $F$ has nontrivial negative normal bundle. By Lemma \ref{xmas2}, $H+\alpha$ restricts nonexactly
 on $H_T(M_{c -\epsilon})$ for $\epsilon >0$ small.

It follows from Tolman-Weitsman, that there is some component of the fixed point set $F' \subset M^T$, for
which $f(F') < c -\epsilon < c$ and $\eta_T|_{F'}$ is nonexact. Iterating this process we find that $\eta_T$ restricts to
 a nonexact form on $M_{min}$. Replacing $f$ with $-f$ gives us the result for $M_{max}$.
\end{proof}


\begin{proof}[Proof of Theorem \ref{th2}]
We let $\mu = f$ and proceed by induction on the regular intervals of $f$.

Let $c_0 < ...< c_n$ denote the critical values of $f$. By Corollary \ref{snoopy}, $H +\alpha$ restricts to something nonexact
 on $M_{min} = M_{c_0}$. Thus by Lemma \ref{xmas2}, it restrict to something nonexact on $M_{c_0+\epsilon}$, and so restricts nonexactly
  to $M_t$ for all $t \in (c_0,c_1)$.

Now suppose that the first interval on which $H+\alpha$ restricts exactly is $(c_i,c_{i+1})$ for some $i>1$. By Lemma \ref{xmas2},
 we deduce that $H+\alpha$ restricts to an exact form on $M_{c_i}^T$. Theorem \ref{thtw} tells us that there exist classes $a,b \in H_T(M)$
  such that $a +b = [H+\alpha]$ and $a$ restricts to zero on $H_T( f^{-1}((-\infty, t))^T)$ and $b$ restricts to zero on
  $H_T(f^{-1}((t,\infty)^T)$ for $t \in (c_i,c_{i+1})$. But since $H+\alpha$ restricts to zero in $H_T(M_{c_i}^T)$, $a, b$
   also work to show that $[H+\alpha]$ restricts to zero on $H_T(f^{-1}(M_t)$ for $t \in (c_{i-1}, c_i)$, which is a contradiction.
   Thus $H+\alpha$ restricts nonexactly to $M_t$ for all regular $t$ (in fact for all $t$).
\end{proof}

\begin{proof}[Proof of Corollary \ref{ld}]
	Torus moment maps can always be shifted by an arbitrary element of $\lie{t}^*$, so we assume without loss of generality that $\theta = 0$.
	
We use induction on the rank of $T$. For $T$ with rank 1, Corollary \ref{ld} holds by Theorem \ref{th2}.
For a rank 1 subtorus, $S \subset T$ with Lie algebra $\lie{s}$, the restricted $S$ action is Hamiltonian with moment map
 $\mu_S = \pi_{\lie{s}} \circ \mu$, where $\pi_{\lie{s}}: \lie{t}^* \rightarrow \lie{s}^*$ is projection.
 We may choose $S$ so that $0 \in \lie{s}$ is a regular value of $\mu_S$ (see Proposition 3.12 in \cite{BL08} for example).
  Thus, $S$ acts freely on $\mu_S^{-1}( 0)$ and the quotient $M //S$ is smooth. Moreover, the induced Hamiltonian action
   of $T/S$ on $M //S$ satisfies the hypotheses of the corollary. So by reduction in stages and induction, the quotient
   $$ M//T \cong (M//S)//(T/S) $$ has nonexact twisting.
\end{proof}

\section{New examples of Hamiltonian generalized complex torus actions}\label{examples}

\subsection{ Logarithmic transformation for generalized complex four manifolds}

An important ingredient of our construction is a surgery introduced by Cavalcanti and Gualtieri \cite{CaGu07}
for $4$-manifolds. This surgery removes a neighborhood of a symplectic $2$-torus
 and replaces it by a neighborhood of a torus where the generalized complex structure changes type, and is an example
  of logarithmic transformation as defined by Gompf and Mrowka \cite{GM93}.

 Let $T \hookrightarrow M$ be a $2$-torus with trivial normal bundle in a $4$-manifold, and let
$U \cong D^2 \times T^ 2$ be a tubular neighborhood, where $D$ denotes the open unit disk in the complex plane $\C$.
Let  $\psi : S^1 \times T^2 \rightarrow \partial U$ be a diffemorphsim from $S^1\times T^2$ to $\partial U$. Then the adjunction space
\begin{equation} \label{log-transform}\widetilde{M} = (M\setminus U) \cup_{\psi}  (D \times T^2) \end{equation}
 is a manifold, and is said to be a $C^{\infty}$ logarithmic transform of $M$.

 Now we are ready to state the following theorem due to Cavalcanti and Gualtieri.

\begin{theorem}(\cite{CaGu07}) \label{surgery-4-GC-mflds} Let $(M, \sigma)$ be a symplectic $4$-manifold,
 $T \rightarrow M$ a symplectic $2$-torus with trivial normal bundle and tubular neighborhood $U$.
 Let $\psi : S^1 \times T^2 \rightarrow \partial U \cong
S^1\times T^2$  be the map given on standard coordinates by
\[\psi(\theta_1, \theta_2, \theta_3) = (\theta_3,  \theta_2, -\theta_1).\]
Then the $C^{\infty}$-logarithmic transform $\widetilde{M}$ of $M$ as defined in (\ref{log-transform})
admits a generalized complex structure $\J$ which exhibits type changing along a $2$-torus, and which
is integrable with respect to a $3$-form $H$, such that $[H]$ is the Poincar\'{e} dual to the
circle in $S^1 \times T^2$ fixed by $\psi$.
\end{theorem}

\begin{remark} When restricting to the open subset $D\times T^2 \subset \widetilde{M}$, the generalized complex
structure $\J$ is exactly the one described in Example \ref{local-model}.
Its type jumps from zero to two along the central fiber $\{0\}\times T^2$.
\end{remark}

\begin{example}(\cite{CaGu07}) \label{four-dim-example} Let $N$ be a two dimensional compact symplectic surface, and $M=N\times T^2$ the product
symplectic manifold, where $T^2$ is the two dimensional torus equipped with the canonical symplectic form.
Fix a point $a \in N$. The image of $T^2$ under
 the embedding \[T^2 \rightarrow V , \,\,\,\,x \mapsto (a, x)\]  is a symplectic two torus with trivial normal bundle in $M$.

Let $\widetilde{M}$ be the adjunction space defined in (\ref{log-transform}). Then by Theorem  \ref{surgery-4-GC-mflds}, $\widetilde{M}$
admits an $H$-twisted generalized complex structure which has type change along a $2$-torus, and which is integrable with respect
to a $3$-form $H$. Moreover, it is easy to see that $\widetilde{M}$ is diffeomorphic to  $X^3\times S^1$ for a compact $3$-manifold $X$.
Note that $X^3$ must be an orientable three manifold,  and that $H$ represents a generator for $H^3(X^3,\R)\cong \R$.
\end{example}

    \subsection{Surgery for higher dimensional generalized complex manifolds}

 Under certain conditions the surgery process described in Example \ref{four-dim-example} can be extended to higher dimensional
 manifolds. Instead of starting with a point in a symplectic two manifold $N$, we are starting with a $2n-4$ dimensional symplectic
 submanifold $A$ of a $2n-2$ dimensional symplectic manifold $N$, which we require to have trivial normal
 bundle. We are going to construct a $2n$ dimensional generalized complex manifold using this initial data.

 By assumption, an open neighborhood $V$ of $A$ in $N$ can be identified with $A\times D$,
  where $D $ is the open unit disk in the complex plane $\C$. Thus $V\times T^2\cong (A\times D)\times T^2$ is
  a tubular neighborhood of $A\times\{0\}\times T^2$; moreover, $\partial (V\times T^2)\cong A\times S^1\times T^2$.

  Let $r,\theta_1$ be the polar coordinates on $D$, $\theta_2,\theta_3$ coordinates for $T^2$ with unit periods,
  and $D_{\frac{1}{\sqrt{e}}}$ the open disk in $\C$ centered at the origin with radius $\dfrac{1}{\sqrt{e}}$.
  Define the map
    \begin{equation}\label{gluing-map} \begin{split} & \psi': A\times (D\setminus D_{\frac{1}{\sqrt{e}}})
    \times T^2  \rightarrow A\times (D\setminus \{0\})\times T^2,
    \\& (x,r, \theta_1,\theta_2,\theta_3) \mapsto (x,\sqrt{\log (er^2)}, \theta_3,\theta_2,-\theta_1),
     \,\,\, x\in A,   \end{split}\end{equation}

    and the adjunction space
    \begin{equation} \label{surgery-GC-mfld} \widetilde{M}= \left( M\setminus V\right) \cup_{\psi'} \left( A\times (D\times T^2)\right).
    \end{equation}

    We are ready to prove the following result by repackaging the argument used in \cite{CaGu07} to prove
    Theorem \ref{surgery-4-GC-mflds}.

   \begin{theorem} \label{higher-dim-surgery} Suppose that $N$ is a $(2n-2)$ dimensional symplectic manifold, $A\subset N$ is a $(2n-4)$ dimensional
   compact symplectic submanifold with trivial normal line bundle, and that $T^2$ is a two dimensional symplectic torus. Let
   $\widetilde{M}$ be the adjunction space defined in (\ref{surgery-GC-mfld}), and $S^1\subset T^2$ the circle parametrized by $\theta_2$.
   Then $\widetilde{M}$ admits an $H$-twisted generalized complex manifold whose type jumps from zero to two along $A\times\{0\}\times T^2$.
   Moreover, the closed three form $H$ represents the Poincar\'{e} dual to $A\times \{0\}\times S^1$.
   \end{theorem}

     \begin{proof} Let $\sigma$ be the symplectic structure on $A$, and $\pi: A\times D\times T^2 \rightarrow A$
     the projection map to the first factor. Then by the symplectic tubular neighborhood theorem,  $V\times T^2$ is symplectomorphic
      to $A\times D\times T^2$ with the symplectic structure
 \[ \omega:=\pi^*\sigma+ \dfrac{1}{2}dr^2\wedge d\theta_1+d\theta_2\wedge d\theta_3.\]
 Now consider the symplectic structure $\alpha:=\pi^*\sigma+ \gamma$ on $A\times (D\setminus D_{\frac{1}{\sqrt{e}}})\times T^2$, where
  $\gamma$ is the symplectic structure on $(D\setminus \{0\}) \times T^2$ defined in (\ref{B-symplectic}).  Then the map
 \[ \psi': (A\times (D\setminus D_{\frac{1}{\sqrt{e}}})
 \times T^2, \alpha) \rightarrow (A\times (D\setminus \{0\})\times T^2, \omega)\] is a symplectomorphism.

  Let $B$ be the closed two form on $(D\setminus \{0\})\times T^2$ defined in (\ref{B-symplectic}). Its pullback
  under the projection map \[A\times \left((D\setminus \{0\})\times T^2\right)\rightarrow (D\setminus \{0\})\times T^2 \] defines a two form
  on $A\times (D\setminus \{0\})\times T^2$, which for simplicity we will still denote by $B$.
  Choose an extension $\tilde{B}$ of $(\psi'^{-1})^*B$ to $M\setminus (A\times \{0\}\times T^2)$. Then the
  pure spinor $e^{\tilde{B}+i\omega}$ determines a generalized complex structure of type zero on
  $M \setminus (A\times \{0\}\times T^2)$, which is integrable with respect to the $d\tilde{B}$-Courant bracket.

 Note that $dB=0$ on $A\times D\times T^2$. So $H:=d\tilde{B}$ is actually a globally defined closed three form
 on $\widetilde{M}$. By Example \ref{local-model}, $D\times T^2$ carries a generalized complex structure whose type jumps
 from $0$ to $2$ along $\{0\}\times T^2$. So $A\times D\times T^2$ carries a generalized complex structure
 whose type jumps from $0$ to $2$ along $A\times \{0\}\times T^2$; moreover, outside the type jumping locus, this generalized
 complex structure is determined by the pure spinor $e^{B+i\alpha}$. However, the gluing map $\psi'$ satisfies $(\psi')^*(\tilde{B}+i\omega)=B+i\alpha$. We conclude that $\widetilde{M}$ admits a
generalized complex structure whose type jumps from $0$ to $2$ along $A\times\{0\}\times T^2$. Furthermore, this structure is
integrable with respect to the $H$-Courant bracket.

The two form $\tilde{B}$ can be chosen so that it vanishes outside a larger tubular neighborhood $A\times D'\times T^2$ of
$A\times \{0\}\times T^2$, and so that $H=d\tilde{B}$ has support in $A\times(D'\setminus D)\times T^2$ and has the form
\[H=f'(r)dr\wedge d\theta_1\wedge d\theta_3,\] for a smooth bump function $f$ such that $f\mid_{A\times D\times T^2}=1$ and
$f=0$ outside $A\times D'\times T^2$.  Therefore $H$ represents the Poincar\'{e} dual to the Cartesian product of $A$ and the circle
parametrized by $\theta_2$. This finishes the proof of Theorem \ref{higher-dim-surgery}.

     \end{proof}

             \subsection{Symplectic toric manifolds}

        It turns out that symplectic toric manifolds provide us many examples which satisfy the assumptions in Theorem \ref{higher-dim-surgery}.
        We collect here some basic facts concerning symplectic toric manifolds which we need for our construction. Let $T^k \cong \R^k/\Z^k$ be the standard rank $k$ torus. We identify $ \lie{t} \cong \R^k$ with standard basis $\{ e^1,...,e^k\}$ and $\lie{t}^* \cong \R^k$ with dual basis $ \{ e_1,...e_k\}$.

  \begin{definition} \label{symplectic-toric} A symplectic toric manifold is a connected symplectic
          manifold $(M, \omega)$ equipped with an effective Hamiltonian action of the standard torus $T^k$
          where $2k$ is the dimension of $M$. We denote the moment map $\mu = (\mu_1,...\mu_k): M \rightarrow \R^k = \lie{t}^*$.
        \end{definition}

     \begin{definition}  \label{equivalent-toric-mflds} Two symplectic toric manifolds, $(M_i, \omega_i, T^k, \mu^i), i = 1, 2$, are
     \textbf{equivalent} if there exists
     symplectomorphism $\varphi : M_1 \rightarrow M_2$ such that $\mu^1 = \mu^2 \circ \varphi$.
     \end{definition}

        \begin{definition}  A Delzant polytope $\bigtriangleup$ in $R^k$ is a polytope satisfying:
        \begin{itemize} \item[1)] simplicity, i.e., there are k edges meeting at each vertex;
                         \item[2)] rationality, i.e., the edges meeting at the vertex p are rational in the sense
                            that each edge is of the form $p + tu_i, t \geq 0$, where $u_i \in Z^n$;
                          \item[3)]  smoothness, i.e., for each vertex, the corresponding $u_1, \cdots, u_n$ can be chosen
                            to be a $Z$-basis of $Z^k$.
         \end{itemize}
        \end{definition}

        \begin{theorem}(\cite{De}) \label{delzant-classification}

        Compact symplectic toric manifolds are classified by Delzant polytopes. More specifically, the bijective correspondence
        between these two sets is given by the moment map:
       \[ \begin{split} &\{\text{ equivalence classes of compact symplectic toric manifolds}\} \xrightarrow{1-1}  \{\text{ Delzant polytopes}\}\\
        &(M^{2n}; \omega, T^k; \mu)\mapsto \mu(M) .\end{split}\]
        \end{theorem}

       \begin{lemma}\label{isotropy-of-toric-mflds} (\cite[Lemma 6.3]{LeTo97}) Let $(M, \omega, T^k, \mu)$ be a compact symplectic
       toric manifold. Then for any $x\in M$, the isotropy group $T_x$  is connected, with Lie algebra $\lie{t}_x = (TF)^{\perp}$, where $\mu(x)$ lies in the interior of the face $F$ of $\bigtriangleup$.  \end{lemma}

        We also need the construction of symplectic cutting which was first introduced by Lerman \cite{Le95}.
        Let $(M, \omega)$ be a symplectic manifold with a Hamiltonian $T^1$ action and a moment map $\mu : M \rightarrow \R$.
         Suppose that $c, d \in \R$ such that $T^1$ acts freely on $\mu^{-1}(c)$ and $\mu^{-1}(d)$. Consider the quotient space
       \[ \overline{M}_{[c,d]} :=\mu^{-1}([c, d])/\sim,    \]
          where $x\sim y$ if and only if $\mu(x)=\mu(y) \in \{c,d\}$ and $y=t \cdot x$ for some element $t \in T^1$.
      The interested readers are referred
       to \cite{Le95} for the proof of the following lemma.
        \begin{lemma} \label{symplectic-cutting-property} $\overline{M}_{[c,d]}$ is a symplectic Hamiltonian $T^1$-space with moment map $\mu_{cut}$ induced by $\mu$, in such a way
       that the obvious homeomorphism $\mu^{-1}((c, d)) \cong \mu^{-1}_{cut}((c,d))$ is a $T^1$-equivariant symplectomorphism intertwining $\mu$ and $\mu_{cut}$. If in addition to the action of a
        circle on $(M,\omega)$ we have a Hamiltonian action
        of another group $K$ on $M$ that commutes with the action of $T^1$, then the space $\overline{M}_{[c,d]}$
        is again a Hamiltonian $K \times T^1$ space.   \end{lemma}

        \begin{definition}\label{symplectic-cut} We call the symplectic manifold $\overline{M}_{[c,d]}$ the (two sided) symplectic cut of $M$ with respect to
        the ray $[c,d]$ (with symplectic form and moment map understood).  \end{definition}

       In fact, Lerman defines the symplectic cut with respect to a ray $[c,\infty)$. The two-sided cut described above is produced by cutting with respect to a ray $[c,\infty)$ and then with respect to $(-\infty,d]$.

    \subsection{Compact examples of Hamiltonian GC manifolds with non-trivial twisting }

\begin{lemma}\label{toric-observation}
	Let $H = span\{e_1,...,e_{k-1}\} \subset R^{k}$ and let $\bigtriangleup$ be a Delzant polytope which intersects $H$ orthogonally in the sense that $H \cap \bigtriangleup$ is nonempty and for every face $F$ of $\bigtriangleup$ that intersects $H$ nontrivally, the tangent space $TF$ contains $e_k$.
	If $(M, \omega, T^k, \mu = (\mu_1,...,\mu_k))$ is the toric manifold associated to $\bigtriangleup$, then for $\delta >0$ sufficiently small there is an equivalence of (noncompact) symplectic toric manifolds $$ \phi: \mu_k^{-1}((-\delta, \delta) ) \cong A \times C$$ where $C \cong (-\delta, \delta) \times S^1$ is a symplectic cylinder with moment map projection onto $(-\delta,\delta)$ and $A$ is the $2k-2$ dimensional toric manifold with moment polytope equal to $H \cap \bigtriangleup$.
\end{lemma}

\begin{proof}
	The hyperplane $H$ does not intersect any vertices of $\bigtriangleup$ because that would contradict the orthogonality condition. Since $\bigtriangleup$ has only a finite number of vertices, it follows that we may choose $\delta >0$ so that the region between $H + \delta e_k$ and $H - \delta e_k$ contains no vertices.
Let $S \subset T^k$ denote circle with Lie algebra the span of $e^k$. The orthogonality condition combined with Lemma \ref{isotropy-of-toric-mflds}
 implies that $S$ acts freely on $\mu^{-1}(H)$ and so, choosing $\delta$ smaller if necessary, $S$ also acts freely on $\mu^{-1}(H \pm \delta e_k) = \mu_k^{-1}( \pm \delta)$. Thus we may form the two-sided symplectic cut with
respect to $S$ to get a compact symplectic toric manifold $ M_{[-\delta, +\delta]}$ with moment polytope the product $ (H\cap \bigtriangleup) \times [-\delta, \delta] $. Delzant's theorem tells us that $M_{[-\delta, +\delta]}$ is equivalent as a Hamiltonian $T^k$ space to $ A \times \CP^1$, where $A$ is as above and $\CP^1$ has symplectic volume $2\delta$ and is acted on by rotation by $S$ about an axis. Finally, $ \mu_k^{-1}((-\delta, \delta) )$ is identified with subset $A \times C \subset A \times \CP^1$, where $C$ equals $\C P^1$ minus the two poles.
\end{proof}

We are ready to state the main result of this section.

\begin{theorem} \label{main-construction} Let
$(N, \omega, T, \Phi)$ be a $2n-2$ dimensional symplectic toric manifold which
satisfies the condition described in Lemma \ref{toric-observation}, where $k = n-1$.
Then\begin{itemize}

\item [(a)]  $N$ contains a compact $(2n-4)$ dimensional symplectic submanifold $A$ with trivial normal bundle.
\item [(b)] Identify a tubular neighborhood of $A$ with $A\times D$, where $D\subset \C$ is the unit
open disk. Then the adjunction space $\widetilde{M}$ as defined in (\ref{surgery-GC-mfld}) admits an $H$-twisted
generalized complex structure $\J$ whose type jumps from $0$ to $2$ along the locus $A\times\{0\}\times T^2$;
moreover, the closed three form $H$ represents the Poincar\'{e} dual of $A\times \{0\}\times S^1$, where
$S^1 \subset T^2$ is the circle parametrized by $\theta_2$ as in (\ref{gluing-map}).

\item [(c)] There exists an $(n-2)$ dimensional torus  which acts on $\widetilde{M}$ in a Hamiltonian fashion.
\end{itemize}
\end{theorem}

\begin{proof} Let $A$ be the $(2n-4)$ dimensional symplectic toric manifold described in Lemma \ref{toric-observation}.
Choose the origin $0 \in C$. By Lemma \ref{toric-observation}, $A\times \{0\}$ embeds onto
a compact $(2n-4)$ dimensional symplectic submanifold of $N$ with trivial normal bundle, which for simplicity
we will still denote by $A$. Assertions (a) and (b) of Theorem \ref{main-construction} follow immediately.

Now let $S$ be the circle described in Lemma \ref{toric-observation}, and $T^{n-2}$
 the complementary subtorus spanned by $\{e^1,...,e^{n-2}\}$ so that $T^{n-1}=T^{n-2}\times S$. Letting $T^{n-2}$ act trivially on $T^2$,
  we extend the Hamiltonian action of $T^{n-2}$ on $N$ to a Hamiltonian action on the symplectic manifold $M=N\times T^2$.
  Observe that this action preserves the tubular neighborhood $A\times \left( D\times T^2\right)$
  of $A\times \{0\}\times T^2$ and is, moreover, Hamiltonian on the first factor $A$
   and trivial on the second factor $D\times T^2$. Now it follows from Example \ref{Hamiltonian-GCmflds} and
   the description of the generalized complex structure $\J$ in the proof of Theorem \ref{higher-dim-surgery}
   that the action of $T^{n-2}$ on $M$ determines naturally a Hamiltonian action on $(\widetilde{M},\J)$. Indeed,
   $\Phi$ extends in a natural way to a moment map $\tilde{\Phi}$ of the generalized complex Hamiltonian space $\widetilde{M}$.

\end{proof}


We complete this section by showing that that the generalized complex Hamiltonian manifolds constructed in Theorem \ref{main-construction} are actually topologically twisted. To do this we first show that the surgery process commutes with reduction.
Let $N$ be a symplectic toric manifold satisfying the hypotheses Theorem \ref{main-construction}, and let $G$ denote a subtorus of $T^{n-2} \subset T^{n-2}\times S= T^{n-1}$. We may consider the restricted Hamiltonian action by $G$ with moment map $\Phi_G: N\rightarrow \frak{g}^*$. If $G$ acts freely on $\phi_G^{-1}(0)$, we may form the generalized complex quotient $N//_{0}G$ which we abbreviate $N//G$ (notice in particular that this holds when $G = T^{n-2}$). Under these hypotheses we have the following.

\begin{proposition}[Surgery commutes with reduction]\label{comm}
The quotient space $N//G$ is a symplectic toric manifold satisfying the conditions specified in Theorem \ref{main-construction}. Let $\widetilde{M//G}$ be the adjunction space as defined in (\ref{surgery-GC-mfld}) with $M$ replaced by $M//G \cong (N//G) \times T^2$.
Then we have an isomorphism of generalized complex Hamiltonian $T^{n-1}/G$-manifolds $$ \widetilde{M//G} \cong \widetilde{M} // G.$$
\end{proposition}

\begin{proof}

Since the surgery is invariant under automorphisms on the $T^{n-2}$ factor, we may assume the $G$ has Lie algebra $Span\{ e^1,..., e^r\}$ for some $r \leq n-2$.  	

Let $\bigtriangleup$ denote the moment polytope of the symplectic toric manifold $N$. The moment polytope for $ N//G$ can be identified with the intersection of the perpendicular of the Lie algebra of $G$ with $\bigtriangleup$, i.e., $ Span\{ e_{r+1},...,e_{n-1}\} \cap \bigtriangleup$, which evidently satisfies the hypotheses of Lemma \ref{toric-observation}. So $N//G$ satisfies the conditions specified in Theorem \ref{comm} and it makes sense to form $\widetilde{M//G}$.

To compare with $\widetilde{M}//G$, notice that the $T^{n-2}$-action, and hence also the $G$-action on $\widetilde{M}$, preserves the partition $ \widetilde{M} = (M \setminus U) \cup_{\psi'}  (A\times D\times T^2)$. After quotienting we obtain a partition $\widetilde{M}//G = (M \setminus U)//G \cup_{\psi'}  (A\times D\times T^2)//G \cong (M//G \setminus U//G) \cup_{\psi'}  ((A//G)\times D\times T^2)$, which is easily identified as the surgery applied to $M//G = N//G \times T^2$.
\end{proof}

\begin{corollary}
The generalized complex manifold $\widetilde{M}$ constructed in Theorem \ref{main-construction} is topologically twisted.
\end{corollary}

\begin{proof}
Apply Proposition \ref{comm} when $G = T^{n-2}$. This is possible because the orthogonality condition on the moment polytope combined with Lemma \ref{isotropy-of-toric-mflds} imply that $G$ acts freely on $\phi_{G}^{-1}(0)$. By the Delzant classification, the toric 2-manifold $N//G$ is diffeomorphic to $\C P^1$. In this case the surgery constructing $\widetilde{M//G}$ the same as Example \ref{four-dim-example}, so it is clear that $\widetilde{M}//G = \widetilde{M//G}$ is topologically twisted. By proposition \ref{strengthening}, this then implies that $\widetilde{M}$ is also topologically twisted.
\end{proof}

\medskip

\noindent
Mathematical Institute, 24-29 St. Giles¡¯, Oxford, OX1 3LB, England,\\ 
Thomas.Baird@maths.ox.ac.uk\\

\noindent

\smallskip

\noindent
Department of Mathematical Sciences, Georgia Southern University, 203 Georgia Ave., Statesboro, GA, 30460, yilin@georgiasouthern.edu

\noindent

\end{document}